\documentclass[12pt]{amsart}

\newtheorem{theorem}{Theorem}[section]
\newtheorem{lemma}[theorem]{Lemma}

\newtheorem{proposition}{Proposition}[section]
\newtheorem{corollary}[theorem]{Corollary}

\theoremstyle{definition}

\theoremstyle{remark}
\newtheorem{remark}[theorem]{Remark}

\numberwithin{equation}{section}

\renewcommand{\a}{\alpha}
\newcommand{\ch}{\mbox{ch}}
\newcommand{\diag}{\mbox{diag}}
\newcommand{\FG}{{\mathcal F}_{\G} }
\newcommand{\FGG}{\overline{\mathcal F}_{\G} }

\newcommand{\G}{\Gamma}
\newcommand{\Gbar}{\overline{\Gamma}^*}
\newcommand{\Gm}{{\Gamma}_m}
\newcommand{\Gn}{{\Gamma}_n}

\newcommand{\gs}{\mathfrak g}
\newcommand{\g}{\gamma}
\newcommand{\hg}{\widehat{\mathfrak h}_{\G, \wt}}
\newcommand{\hh}{\widehat{\mathfrak h}}
\newcommand{\hhg}{\widehat{\widehat{\mathfrak g}}}
\newcommand{\la}{\lambda}
\newcommand{\loopg}{\widehat{\mathfrak g}}
\newcommand{\LG}{ {\Lambda}_{\G}}

\newcommand{\ep}{\epsilon}
\newcommand{\RG}{R_{ \G}}
\newcommand{\RGG}{\overline{R}_{ \G}}
\newcommand{\RGO}{R^0_{ \G}}
\newcommand{\Rz}{R_{\mathbb Z} (\G)}
\newcommand{\Rzz}{\overline{R}_{\mathbb Z} (\G)}
\newcommand{\SG}{ S_\G }
\newcommand{\SGG}{\overline{S}_\G }
\newcommand{\SGO}{ S^0_\G }

\newcommand{\VG}{V_{ \G}}
\newcommand{\VGG}{\overline{V}_{ \G}}

\newcommand{\wt}{\xi}
\newcommand{\Z}{\mathbb Z}

\begin{document}

\title[Vertex representations and McKay correspondence]
    {Vertex representations via finite groups and
the McKay correspondence
    }
\author{Igor B. Frenkel}
\address{Frenkel: Department of Mathematics,
   Yale University,
   New Haven, CT 06520}
\author{Naihuan Jing}
\address{Jing: Department of Mathematics,
   North Carolina State Univer\-sity,
   Ra\-leigh, NC 27695-8205}
\email{jing@math.ncsu.edu}
\thanks{Research of I.B.F. is supported by
NSF grant DMS-9700765;
research of N.J. is supported in part by NSA grant
MDA904-97-1-0062.}
\author{Weiqiang Wang}
\address{
Wang: Department of Mathematics,
   Yale University,
   New Haven, CT 06520}
\email{wqwang@math.yale.edu}

\maketitle

\section{Introduction} \label{S:intro}
    It was first realized in \cite{FK} and \cite{S1} that
the basic representation $V$ of an affine Lie algebra $\loopg$
of ADE type can be constructed from the root lattice $Q$ of
the corresponding finite dimensional Lie algebra $\gs$  as follows:

\begin{eqnarray}   \label{eq_space}
  V = S (Q \bigotimes \mathbb C [ t^{ -1}] t^{ -1})
      \bigotimes \mathbb C [Q] ,
\end{eqnarray}
where the first factor is a symmetric algebra and the second
one is a group algebra.
The affine algebra $\loopg$ contains a Heisenberg algebra
$\hh$. One can define the so-called vertex operators
$X(\alpha, z)$ associated to $ \alpha \in Q$ acting on $V$
essentially using the Heisenberg algebra $\hh$.
The representation of $\loopg$ on $V$ is
then obtained from the action of the Heisenberg algebra $\hh$ and
the vertex operators $X(\alpha, z)$ associated to $\alpha$
in the root system of $\gs$.

This construction was extended in \cite{F2} to more general
lattices to provide vertex representations of affinization
of Kac-Moody Lie algebras. The special case of affinization
of affine Lie algebras called toroidal Lie algebras was
discussed further in \cite{MRY}.

Another important special case of vertex representations
is given by the standard lattice $\mathbb Z^N$. The vertex operators
corresponding to the unit vectors in $\mathbb Z^N$ give
rise to a representation of an infinite-dimensional Clifford
algebra \cite{F1}, the relation known as boson-fermion
correspondence. In the special case $N =1$, the
transition matrix between the monomial bases for
representations of Heisenberg and Clifford algebras
yields the character tables of symmetric groups $S_n$ for all $n$
\cite{F3} (see \cite{J1}).

Vertex operators have appeared and
played an important role in many diversed
fields of mathematics and physics (cf. \cite{FLM} and references
therein).  Surprisingly, they also appeared in recent works of
Nakajima (cf. \cite{N} and references therein)
and independently of Grojnowski \cite{Gr}
on the (homology groups of) Hilbert schemes of points on a surface $X$.
In particular when $X$ is a minimal resolution
$\widehat{\mathbb C^2 /\G}$ of
the simple quotient singularity of $\mathbb C^2$ by a finite subgroup
$\G$ of $SU_2$, they were able to realize geometrically
the vertex representations of affine Lie algebras.

Vertex operators showed up in equivariant K-theory
as well. Partly motivated by a footnote in \cite{Gr},
Segal \cite{S2} outlined
the $S_n$-equivariant K-theory of the $n$-th direct product $X^n$ of $X$.
Generalizing \cite{S2}, the third-named author \cite{W}
studied in detail the $\Gn$-equivariant K-theory of $X^n$
for a $\G$-space $X$, where $\Gn$ is the so-called
wreath product which is the semi-direct product
of the symmetric group $S_n$ and
the $n$-th direct product of a finite group $\G$.
Among other results a construction of
vertex representations was proposed in \cite{W}
in terms of representations of
wreath products; in particular when $\G$ is a subgroup of $SU_2$,
a realization of affine algebras thus obtained can be regarded
as a new form of McKay correspondence \cite{Mc}.

The goal of this paper is to establish firmly such
a finite group theoretic approach
of the vertex representations corresponding to
a very general class of lattices, including the (affine) root
lattices of ADE type and the standard lattice $\mathbb Z^N$
as special cases.
It was conjectured \cite{W} that there exists a natural isomorphism from the 
representation ring of $\Gn$ for
a finite subgroup $\G$ of $SU_2$ to the homology group
of the Hilbert scheme of $n$ points on $\widehat{\mathbb C^2 /\G}$
after the dimensions of these two spaces
were shown to be the same. In this way,
the group theoretic realization of vertex representations
for affine algebras in this paper might be regarded as a counterpart of
the geometric construction given in \cite{Gr, N} for the surface
$\widehat{\mathbb C^2 /\G}$.

Let us explain the contents of this paper in more detail.
Let $\G$ be an arbitrary finite group. Denote by $R_{\Z} (\Gn)$
the representation ring of the wreath product $\Gn$ and
let $R(\Gn) = R_{\Z} (\Gn) \otimes_{\mathbb Z} \mathbb C$. We set
\begin{eqnarray*}
  \RG = \bigoplus_{n \geq 0} R (\Gn).
\end{eqnarray*}
The representation theory of wreath products
was developed by Specht (cf. \cite{M2, Z}).
Given a self-dual virtual character $\wt$ in $R(\G)$,
we introduce
a {\em weighted} bilinear form on $R(\G_n)$
and an induced one on $\RG$ denoted by $\langle \ , \ \rangle_{\wt}$.
When $\wt$ is trivial, it reduces to the standard
bilinear form on the representation ring of
a finite group.

It turns out that the spaces constructed from the representation
theory of $\Gn$ reproduce those appearing in
vertex representations (see Eqn.~(\ref{eq_space})). In fact one
can construct the symmetric algebra $\SG$ associated to
the lattice $\Rz$. We define a natural bilinear
form $\langle \ , \ \rangle_{\wt}'$ induced from
the $\wt$-weighted bilinear form on $\RG$. We define
a characteristic map $\ch$ from
$(\RG, \langle \ , \ \rangle_{\wt}) $ to
$(\SG, \langle \ , \ \rangle_{\wt}') $ and show
that it is an isometry. The map $\ch$
for $\wt$ being a trivial character was considered in \cite{M1, M2}.

We introduce an infinite-dimensional Heisenberg
algebra $\hg$ acting on $\SG$ associated to $\G$ and $\wt$, and thus
regard $\SG$ as the Fock space of $\hg$. We give a group
theoretic realization
\cite{W} of this Heisenberg algebra acting
on $\RG$ which is compatible with the one
on $\SG$ via the characteristic map $\ch$.

Denote by $\FG$ the tensor product of $\RG$ and
the group algebra of $\Rz$ which is
the full space of a vertex representation.
Given $\g \in \Rz $,
we define a generating function $X(\g, z)$ of
group theoretic operators
acting on $\RG$. We prove that $X(\g, z)$ thus defined
is precisely the vertex operator associated to
$\g $ which can be essentially constructed in terms of
the Heisenberg algebra $\hg$.
We then calculate the product of two
vertex operators of the form $X( \g, z), \g \in \Rz$, which
are equivalent to the commutation relations among the
components of $X( \g, z)$.

Now we specialize to the important case when $\G$ is
a finite group of $SU_2$ and $\wt$ is twice the trivial character
minus the character of a two-dimensional natural
representation of $\G$.
According to McKay \cite{Mc}, the lattice
$(\Rz, \langle \ , \ \rangle_{\wt })$ is
positive semi-definite and its radical
is one-dimensional with a generator $\delta$ given by the
regular character of $\G$. Furthermore
$(\Rz, \langle \ , \ \rangle_{\wt })$ can be
identified with the root lattice of the
affine Lie algebra $\loopg$ associated to
a simple Lie algebra $\mathfrak g$ of ADE type.
Our general construction when specialized to
this case gives a group theoretic realization
of the toroidal Lie algebra associated to $\mathfrak g$.

The sublattice $\Rzz$ of
$\Rz$ generated by non-trivial irreducible characters
of $\G$ can be identified with the root lattice of $\mathfrak g$.
It follows that the lattice
$\Rz$ is a direct sum of $\Rzz$ and the lattice $\mathbb Z \delta$
with zero bilinear form generated by $\delta$.
Therefore $\FG$ can be decomposed as a tensor product
of the space $\FGG$ associated to the lattice $\Rzz$ with the space
associated to the lattice $\mathbb Z \delta$.
This gives a new form of ``McKay correspondence'':
starting with a finite subgroup $\G \subset SU_2$ we
have been able to construct the basic representation on $\FGG$
of the affine Lie algebra $\loopg$ with the
Dynkin diagram corresponding to $\G$.

If we specialize $\wt$ to be the trivial character of $\G$,
then the lattice $( \Rz, \langle \ , \ \rangle_{\wt})$ is
just the standard lattice $\mathbb Z^N$ where $N$ is
the number of irreducible characters of $\G$.
The vertex operators corresponding to irreducible
characters of $\G$ generate an infinite-dimensional
Clifford algebra. In this case the transition
matrix between the monomial bases of Heisenberg and
Clifford algebras yields the character table of the
wreath product $\Gn$ for all $n$, generalizing
the symmetric group picture.

We announce some further directions we are currently pursuing.
One is on the $q$-deformation
of our construction while another is on the
generalized McKay correspondence and
affine Lie algebras of non-simply-laced type.
We will also generalize the results in \cite{J1, J2} to
the wreath product setting. In the end we pose
the very important problem
to give a group theoretic construction of the whole
vertex algebra structure on $\FG$ \cite{B, FLM}.

The paper is organized as follows. In Sect.~\ref{sect_wreath}
we review the theory of wreath products.
In Sect.~\ref{sect_weight} we introduce
the weighted bilinear form on $\RG$.
In Sect.~\ref{sect_heis} we define the Heisenberg algebra
$\hg$ and its Fock space together with some natural
basis. In Sect.~\ref{sect_isom} we establish
the isometry between $\RG$ and $\SG$.
In Sect.~\ref{sect_vertex} we define the
vertex operators acting on $\RG$.
In Sect.~\ref{sect_ade} we calculate
products of vertex operators and
establish a new McKay correspondence.
In Sect.~\ref{sect_char} we derive the
character tables of $\Gn$ from the vertex operator
approach.
\section{Preliminaries on wreath products} \label{sect_wreath}
\subsection{The wreath product $\Gn$}
  Let $\Gamma$ be a finite group with $r +1$ conjugacy classes.
We denote by $\Gamma^*$ the set of complex
irreducible characters and by $\Gamma_*$ the set of
conjugacy classes. The character value $ \gamma(c) $
of $\g \in \G^*$ at a conjugacy class $c\in \Gamma_*$
yields the character table $ \{ \gamma(c) \} $ of $\G$.

The space of class functions on $\Gamma$ is given by
$$
R(\Gamma)=\bigoplus_{i= 0}^r \mathbb C \gamma_i,
$$
where $\gamma_i$ are all the inequivalent irreducible characters of $\G$
and in particular $\g_0$ is the trivial character.
We denote by $\Rz$ the integral combination of
irreducible characters of $ \G$. Denote by $ c^i, i = 0, \ldots , r$
the distinct conjugacy classes in $\G_*$, in particular $c^0$ is
the identity conjugacy class.

For $c \in \G_*$
let $\zeta_c$ be the order of the centralizer of an element
in the class $c$, so the
order of the class is then $|c|=|\G |/\zeta_c$.
The usual bilinear form on $R(\G )$ is defined
as follows:

\begin{eqnarray*}
\langle f, g \rangle_{\G} = \frac1{ | \G |}\sum_{x \in \Gamma}
          f(x) g(x^{ -1})
 = \sum_{c \in \Gamma_*} \zeta_c^{ -1} f(c) g(c^{ -1}),
\end{eqnarray*}
where $c^{ -1}$ denotes the conjugacy class
$\{ x^{ -1}, x \in c \}$. Clearly $\zeta_c = \zeta_{c^{-1}}$.
We will refer to this bilinear form as the {\em standard} one
on $R( \G)$. We will often write $\langle \ , \ \rangle$
for $\langle \ , \ \rangle_{\G }$ when no ambiguity may arise.
It is well known that
\begin{eqnarray}
  \langle \g_i, \g_j \rangle &= & \delta_{ij} , \nonumber \\
  \sum_{ \g \in \G^*} \g (c ')  \g ( c^{ -1})
    &= & \delta_{c, c '} \zeta_c, \quad c, c ' \in \G_*.  \label{eq_orth}
\end{eqnarray}
Thus $\Rz$ endowed with this bilinear form
becomes an integral lattice in $R( \G )$.

\begin{remark}
Besides the standard bilinear form on a finite group
one can introduce a standard Hermitian form which
is also useful in our constructions.
Throughout the paper we have chosen to use bilinear
forms instead of Hermitian forms.
\end{remark}

Given a positive integer $n$, let $\Gamma^n = \Gamma \times \cdots
\times \Gamma$ be the $n$-th direct product of $\Gamma$, and let
$\G^0 $ be the trivial group. The symmetric group
$S_n$ acts on $\Gamma^n$ by permutations:
$\sigma (g_1, \cdots, g_n)
  = (g_{\sigma^{ -1} (1)}, \cdots, g_{\sigma^{ -1} (n)}).
$
The wreath product of $\Gamma$ with $S_n$ is defined to be
the semi-direct product
$$
 \Gamma_n = \{(g, \sigma) | g=(g_1, \cdots, g_n)\in {\Gamma}^n,
\sigma\in S_n \}
$$
 with the multiplication
$$
(g, \sigma)\cdot (h, \tau)=(g \, {\sigma} (h), \sigma \tau ) .
$$
Note that $\Gn = S_n$ when $\G =1$, and $\Gn $ is the hyperoctahedral group
of rank $n$ when $\G = \mathbb Z / 2 \mathbb Z$.
\subsection{Conjugacy classes of $\Gn$}
Let $\la=(\la_1, \la_2, \cdots, \la_l)$ be a partition
of integer $|\la|=\la_1+\cdots+\la_l$, where
$\la_1\geq \dots \geq \la_l \geq 1$.
The integer $l$ is called the {\em length} of the partition
$\la $ and is denoted by $l (\la )$.
We will identify the partition $(\la_1, \la_2, \cdots, \la_l)$ with
$(\la_1, \la_2, \cdots, \la_l, 0, \cdots, 0)$.
We will also make use of another notation for partitions:
$$
\la=(1^{m_1}2^{m_2}\cdots) ,
$$
where $m_i$ is the number of parts in $\la$ equal to $i$.

We will use partitions indexed by $\G_*$ and $\G^*$. For
a finite set $X$ and $\rho=(\rho(x))_{x\in X}$ a family
of partitions indexed by $X$, we write
$$\|\rho\|=\sum_{x\in X}|\rho(x)|.$$
Sometimes it is convenient to regard $\rho=(\rho(x))_{x\in X}$
as a partition-valued function on $X$.
We denote by $\mathcal{P}(X)$ the set of all partitions indexed by $X$
and by $\mathcal{P}_n(X)$ the set of all partitions
in $\mathcal{P}(X)$ such that $\|\rho\|=n$.

The conjugacy classes of ${\Gamma}_n$ can be described in the following way.
Let $x=(g, \sigma )\in {\Gamma}_n$, where
$g=(g_1, \cdots, g_n) \in {\Gamma}^n,$ $ \sigma \in S_n$.
The permutation $\sigma $
is written as a product of disjoint cycles. For each
such cycle $y=(i_1 i_2 \cdots i_k)$
the element $g_{i_k} g_{i_{k -1}}
\cdots g_{i_1} \in \Gamma$ is determined up to conjugacy
in $\Gamma$ by $g$ and $y$, and will be called the {\em cycle-product}
of $x$ corresponding to the cycle $y$. For any
conjugacy class $c$ and each integer $i\geq 1$,
the number of $i$-cycles in $\sigma$ whose cycle-product lies in $c$
will be denoted by $m_i(c)$. Denote by $\rho (c)$ the partition
$(1^{m_1 (c)} 2^{m_2 (c)} \ldots )$, $c \in \G_*$.
Then each element
$x=(g, \sigma)\in {\Gamma}_n$ gives rise to a partition-valued
function $( \rho (c))_{c \in \G_*} \in {\mathcal P} ( \G_*)$ such that
$\sum_{i, c} i m_i(c) =n$. The partition-valued function
$\rho =( \rho(c))_{ c \in G_*} $
is called the {\em type} of $x$. It is known (cf. \cite{M2})
that any two elements of ${\Gamma}_n$ are conjugate
in ${\Gamma}_n$ if and only if they have the same type.

Given a partition $\lambda = (1^{m_1} 2^{m_2} \ldots )$,
we define
\[
  z_{\la } = \prod_{i\geq 1}i^{m_i}m_i!.
\]
We note that $z_{\la }$ is the order of the centralizer
of an element of cycle-type $\la $ in $S_{|\la |}$.
The order of the centralizer of an element
$x = (g, \sigma) \in {\Gamma}_n$ of the type
$\rho=( \rho(c))_{ c \in \G_*}$ is
$$
Z_{\rho}=\prod_{c\in \G_*}z_{\rho(c)}\zeta_c^{l(\rho(c))}.
$$
\subsection{Hopf algebra structure on $\RG$}
      We define
\[
  \RG = \bigoplus_{n\geq 0} R({\Gamma}_n).
\]
 $\RG$ carries a natural Hopf algebra structure (cf. e.g. \cite{M1, Z, W})
with multiplication defined by the composition
\[
 m: R(\Gn ) \bigotimes R(\Gm)
 \stackrel{\cong }{\longrightarrow} R(\Gn \times \Gm)
 \stackrel{Ind}{\longrightarrow} R( {\Gamma}_{n + m}),
\]
and comultiplication defined by the composition
\[
 \Delta: R(\Gn ) \stackrel{Res}{\longrightarrow}
 \bigoplus_{m =0}^n R( {\Gamma}_{n - m} \times \Gm)
 \stackrel{\cong }{\longrightarrow}
 \bigoplus_{m =0}^n R( {\Gamma}_{n - m}) \bigotimes R(\Gm) .
\]
Here $Ind: R(\Gn \times \Gm ) \longrightarrow R(\G_{n +m})$
is the induction functor
and $Res: R(\Gn ) \longrightarrow
 R( {\G }_{n - m} \times \Gm )$ is the restriction functor.

The standard bilinear form
in $\RG$ is defined in terms of those on $R( \Gn )$
as follows:
\[
\langle u, v \rangle
 = \sum_{ n \geq 0} \langle u_n, v_n \rangle_{\Gn},
\]
where
$u = \sum_n u_n$ and $v = \sum_n v_n$ with $u_n, v_n\in \Gn$.
\section{Weighted bilinear forms on $R(\G)$ and $R(\Gn)$}
\label{sect_weight}
\subsection{A weighted bilinear form on $R(\G)$}
Let us fix a class function $ \wt \in R(\G)$.
The multiplication in $ R(\G)$ corresponding to the
tensor product of two representations will be denoted by $* $.

We denote by $a_{ij} \in \mathbb C$ the (virtual) multiplicities
of $\g_j$ in $ \wt * \g_i $. In other words, we have
the following decomposition
\begin{eqnarray}  \label{eq_tens}
 \wt  *  \g_i
  = \sum_{j =0}^r a_{ij} \g_j.
\end{eqnarray}
We denote by $A$ the $ (r +1) \times (r +1)$ matrix
$ ( a_{ij})_{0 \leq i,j \leq r}$.

We introduce the following weighted bilinear form
$$
  \langle f, g \rangle_{\wt } = \langle \wt * f ,  g \rangle_{\G },
   \quad f, g \in R( \G).
$$
We also have an alternative formula:
\begin{eqnarray}
  \langle f, g \rangle_{\wt }
   & =& \frac 1{ |\G|} \sum_{ x \in \G} \wt (x)f(x) g (x^{ -1})
     \nonumber  \\
   & =& \sum_{c \in \G_*} \zeta_c^{ -1} \wt (c) f(c) g (c^{ -1}).
     \label{eq_twist}
\end{eqnarray}
In other words, $\langle \ , \ \rangle_{\G}$ is the average of
the character $ \wt * f * \overline{g}$.

In particular, Eqn.~(\ref{eq_tens}) can be reformulated as
$$
   \langle \g_i, \g_j \rangle_{\wt } = a_{ij}.
$$

Throughout this paper we will always
assume that $\wt $ is a {\em self-dual},
i.e. $\overline{\wt } = \wt  $, or equivalently
$\wt (x) = \wt (x^{-1}), x \in \G$.
The self-duality of $\wt$ together with (\ref{eq_twist})
implies that
$$
 a_{ij} = a_{ji},
$$
i.e. $A$ is a symmetric matrix. An equivalent formula
for $\langle \ , \ \rangle_{\wt}$ is

\begin{eqnarray}
  \langle f, g \rangle_{\wt }
   & =& \sum_{c \in \G_*} \zeta_c^{ -1} \wt (c) f(c^{-1}) g (c).
     \label{eq_real}
\end{eqnarray}

\begin{remark}
 If $\wt$ is the trivial character $\g_0$, then the
 weighted bilinear form becomes the standard one on $ R( \G)$.
\end{remark}
\subsection{A weighted bilinear form on $R( \Gn)$}

Given a representation $V$ of $\G$ with character
$\g \in R(\G)$, the $n$-th outer tensor product
$V^{ \otimes  n} $ of $V$ can be regarded naturally as
a representation of the wreath product $\Gn$
whose character will be denoted by $\eta_n ( \g )$:
the direct product $\G^n$ acts on $\g^{\otimes n}$
factor by factor while $S_n$ by permuting the $n$ factors.
Denote by $\varepsilon_n$ the (1-dimensional) sign representation
of $\Gn$ on which $\G^n$ acts trivially while
$S_n$ acts as sign representation.
We denote by $\varepsilon_n ( \g ) \in R(\Gn)$
the character of the tensor product
of $\varepsilon_n$ and $V^{\otimes n}$.

We may extend naturally $\eta_n$ to
a map from $R(\G)$ to $R(\Gn)$ (cf. \cite{W}).
In particular, if $\beta$ and $\g $ are characters of
representations $V$ and $W$ of $\G$ respectively, then
\begin{eqnarray}  \label{eq_virt}
  \eta_n (\beta - \g) =
  \sum_{m =0}^n ( -1)^m Ind_{\G_{n -m} \times \Gm }^{\Gn}
   [ \eta_{n -m} (\beta) \otimes \varepsilon_m (\g ) ] .
\end{eqnarray}

We introduce a {\em weighted bilinear form} on $R( \Gn)$ by letting
$$
  \langle  f, g\rangle_{\wt, \Gn } =
   \langle \eta_n (\wt ) * f, g \rangle_{\Gn} ,
   \quad f, g \in R( \Gn).
$$
We shall see in Corollary~\ref{cor_char} that
$\eta_n (\wt)$ is self-dual. It follows that the
bilinear form $\langle \ , \ \rangle_{\wt}$ is symmetric.

\begin{remark}
 When $n =1$, this weighted bilinear form obviously reduces
 to the weighted bilinear form defined on $R( \G)$.
\end{remark}

On $\RG = \bigoplus_{n} R(\Gn)$ a symmetric
bilinear form is given by
\[
\langle u, v \rangle_{\wt}
 = \sum_{ n \geq 0} \langle u_n, v_n \rangle_{\wt, \Gn } ,
\]
where
$u = \sum_n u_n$ and $v = \sum_n v_n$ with $u_n, v_n\in \Gn$.
\subsection{A McKay specialization}  \label{subsec_mcka}
    Denote by $d_i = \g_i (c^0)$ the dimension of
the irreducible representation of $\G$
corresponding to the character $\g_i$. The following
proposition appeared in \cite{St} which was motivated by
a special important case observed first by McKay \cite{Mc}.

\begin{proposition}
    The column vector
 $$
   v_i = ( \g_0 (c^i), \g_1 (c^i), \ldots, \g_r (c^i ) )^t
  \quad ( i =0 , \ldots, r)
 $$
 is an eigenvector of the matrix $A$
 with eigenvalue $ \wt (c^i )$.
 In particular $(d_0, d_1, \ldots, d_r)^t $ is an eigenvector
 of $A$ with eigenvalue $\wt(c^0)$.
\end{proposition}

Denote by $E$ the character matrix $ [ v_0, v_1, \ldots , v_r ]$
and $D$ the diagonal matrix
$\diag (\wt (c^0), \ldots, \wt (c^r))$.
The proposition above can be reformulated as
$$
 A E = E D.
$$

  We fix an irreducible faithful complex
representation $\pi$ of $\G$ of
dimension $d$. We still assume that $\pi$ is self-dual.
We set
$$
  \wt = d \g_0 - \pi.
$$

In this case, the weighted bilinear form
$\langle \ , \ \rangle_{\wt}$ on $\RG$
become semi-positive definite. The radical of this
bilinear form is one-dimensional and spanned by
the character of the regular representation of $\G$
$$
  \delta = \sum_{i = 0}^r d_i \g_i.
$$

We further specialize to the case in which $\pi$ embeds
$\G$ into $SU_2$, i.e. $d =2$. The classification of finite subgroups
of $SU_2$ is well known. The following is a complete list of them:
the cyclic, binary dihedral, tetrahedral, octahedral
and icosahedral groups.
We recall the following well-known facts \cite{Mc}:
$a_{ii} =2$ for all $i$;
if $ \G \neq \mathbb Z / 2 \mathbb Z$
and $ i \neq j$ then $a_{ij} = 0 $ or $-1$. If
$ \G \equiv  \mathbb Z / 2 \mathbb Z$ then $a_{01} = -2$.

We associate a diagram with vertices
corresponding to elements $\g_i$ in $\G^*$. We draw one edge
(resp. two edges) between
the $i$-th and $ j$-th vertices if $a_{ij} = -1$ (resp. $-2$).
According to McKay \cite{Mc}, the associated diagram can be identified
with affine Dynkin diagram of ADE type and the matrix $A$ is the
corresponding affine Cartan matrix.
\section{Heisenberg algebras and $\Gn$} \label{sect_heis}
\subsection{Heisenberg algebra $\hg $}
      Let $\hg $ be the infinite dimensional
Heisenberg algebra over $\mathbb C$, associated with
$\Gamma$ and $\wt$,
with generators $a_m(\gamma), m \in \Z, \gamma \in\G^*$
and a central element $C$ subject to the following commutation relations:
\begin{equation}  \label{eq_heis}
[a_m( \gamma), a_n(\gamma ')]
 = m \delta_{m, -n} \langle \gamma, \gamma' \rangle_{\wt } C,
 \quad m, n \in \Z, \, \g, \g ' \in \G^*.
\end{equation}
It is convenient to extend $a_m (\g )$ to all
$\g = \sum_{ i =0}^r s_i \g_i \in R(\G )$ $(s_i \in \mathbb C)$
by linearity:
$ a_m ( \g ) = \sum_i s_i \, a_m (\g_i )$.

Denote by $R_0$ the radical of the bilinear form
$\langle \cdot , \cdot \rangle_{\wt }$ in $R(\G)$.
We note that the Heisenberg algebra may contain a large center
since the bilinear form $\langle \cdot , \cdot \rangle_{\wt }$
may degenerate. In this sense we have abused the notion
of Heisenberg algebra.
The center of $\hg $ is spanned by $C$ together with
$a_m ( \g), \g \in R_0, m \in \mathbb Z$.

For $m\in\mathbb Z, c \in \G_*$ we define
$$
 a_{ m}( c) = \sum_{ \g\in \G^*} \gamma(c^{-1}) a_m( \g ).
$$

 From the orthogonality of the irreducible
characters of $\G $ (\ref{eq_orth}) it follows that
\begin{eqnarray*}
  a_m( \g )
   = \sum_{c \in \G_*} \zeta_c^{ -1}
   \g  (c) a_m(c).
\end{eqnarray*}

 Thus $ a_{n}(c)$ ($ n\in \mathbb Z, c \in \G_*$)
and $C$ form
a new basis for the Heisenberg algebra $\hg$.

\begin{proposition}  \label{prop_orth}
 The commutation relations among the new basis for $\hg$
 are given by
 \begin{eqnarray*}
  [ a_m( {c'}^{ -1}), a_n( c )]
    & =& m \delta_{m, -n}\delta_{c', c} \zeta_c \wt (c) C,
   \quad c, c' \in \G_*.
\end{eqnarray*}
\end{proposition}

\begin{proof} We calculate by using Eqns. (\ref{eq_heis}),
 (\ref{eq_real}) and (\ref{eq_orth}) that
 \begin{eqnarray*}
  && [ a_m( {c'}^{ -1}), a_n( c)] \\
  & =& \sum_{\g', \g \in \G^*} \g '(c') \g  (c^{-1})
       [ a_m( \g '), a_n ( \g )]          \\
  & =& m \delta_{m, -n}
         \sum_{\g', \g \in \G^*} \g '(c') \g (c^{-1})
         \langle \g ', \g \rangle_{\wt} C    \\
  & =& m \delta_{m, -n}
         \sum_{\g', \g \in \G^*} \g '(c') \g (c^{-1})
         \biggl (
          \sum_{ \nu \in \G_*} \zeta_{\nu}^{ -1} \wt (\nu)
          \g (\nu) \g '  (\nu^{-1})
         \biggr ) C    \\
  & =& m \delta_{m, -n}
         \sum_{ \nu \in \G_*} \zeta_{\nu}^{ -1} \wt (\nu)
         \biggl (
           \sum_{\g' \in \G^*} \g '(c') \g ' (\nu^{-1})
         \biggr )
         \biggl (
           \sum_{\g \in \G^*} \g (c^{-1})  \g (\nu)
         \biggr )  C    \\
  & =& m \delta_{m, -n}
         \sum_{ \nu \in \G_*} \zeta_{\nu}^{-1} \wt (\nu)
         \delta_{c', \nu} \zeta_{\nu}
         \delta_{c, \nu} \zeta_{\nu} C    \\
  & =& m \delta_{m, -n} \delta_{c, c'} \zeta_c \wt (c) C.
\end{eqnarray*}
\end{proof}
\subsection{Action of $\hg$ on $\SG$ and $\SGG$}
     Denote by $\SG $ the symmetric algebra
generated by $a_{-n}(\g), n \in \mathbb N, \g\in \Gamma^*$.
There is a natural degree operator on $\SG $
$$
  \deg (a_{ -n}( \g)) = n ,
$$
which makes $\SG $ into a $\mathbb Z_+$-graded space.

We define an action of $\hg$ on $\SG$ as follows:
$a_{-n}( \g), n >0$ acts as
multiplication operator on $ \SG $ and $C$ as the identity
operator; $a_n (\g),$ $ n \geq 0$ acts as a derivation of algebra
\begin{eqnarray*}
  a_n (\g). a_{-n_1}( \alpha_1) a_{-n_2} (\alpha_2)
    \ldots a_{-n_k}( \alpha_k)   \\
 = \sum_{i =1}^k \delta_{n,n_i}
 \langle \g , \alpha_i \rangle_{\wt }
   a_{-n_1}( \alpha_1) a_{-n_2}(\alpha_2) \ldots
  \check{a}_{-n_i}( \alpha_i) \ldots a_{-n_k}(\alpha_k )  .
\end{eqnarray*}
Here $n_i > 0, \alpha_i \in R(\G)$ for $i =1, \ldots , k$,
and $\check{a}_{-n_i}( \alpha_i)$ means the very term
is deleted. In other word, the operator
$a_n (\g ), n > 0, \g \in R^0$ acts as $0$, and
$a_n (\g ), n > 0, \g \in R(\G) - R^0$ acts as certain non-zero
differential operator. Note that $\SG$ is not an irreducible
representation over $\hg$ in general since the bilinear
form $\langle \ , \ \rangle_{\wt}$ may be degenerate.

Denote by $\SGO$ the ideal in the symmetric algebra $ \SG$
generated by $a_{-n}(\g), n \in \mathbb N, \g \in R_0$.
Denote by $\SGG$ the quotient $\SG /\SGO$.
It follows from the definition that $\SGO$ is a subrepresentation
of $\SG$ over the Heisenberg algebra $\hg$.
In particular, this induces a Heisenberg algebra action on $\SGG$
which is irreducible.

We denote by $1$ the unit in the
symmetric algebra $\SG$. By abuse of notation we will
denote by $1$ its image in the quotient $\SGG$.
The element $1$ is the highest weight vector in $\SGG$.
\subsection{The bilinear form on $\SG $}
       The space $ \SG $ admits a bilinear form
$\langle \ ,  \ \rangle_{\wt } '$ characterized by
\begin{eqnarray}  \label{eq_bili}
 \langle 1, 1 \rangle_{\wt}' = 1, \quad
 a_n(\g)^* = a_{-n}(\g), \qquad n\in \mathbb Z \backslash 0.
\end{eqnarray}
Here $a_n(\g)^*$ denotes the adjoint of $a_n(\g)$.

  For any partition $\la =( \la_1, \la_2, \dots)$ and $\g \in \G^*$,
we define
$$
  a_{-\la}( \g) = a_{-\la_1}( \g)a_{ - \la_2}( \g) \dots .
$$
For $\rho = ( \rho (\g) )_{ \g \in \G^*} \in {\mathcal P}(\G^* )$,
we define
$$
  a_{ - \rho } = \prod_{\g \in \G^*}  a_{ - \rho (\g)}(\g).
$$
It is clear that $a_{ - \rho},
\rho \in {\mathcal P}(\G^* )$ consist of a
$\mathbb C$-basis for $\SG $.

Given a partition $ \la = ( \la_1, \la_2, \ldots )$
and $c \in \G_*$, we define
\begin{eqnarray*}
  a_{ - \la } (c )  = a_{ - \la_1}(c) a_{ - \la_2 } (c) \ldots.
\end{eqnarray*}
For any $\rho = ( \rho (c) )_{ c \in \G_* } \in
 \mathcal P ( \G_* )$, we further define
\begin{eqnarray*}
  a_{- \rho}' = \prod_{ c \in \G_*} a_{ - \rho (c)} (c).
\end{eqnarray*}
The elements $a_{ - \rho} ',
\rho \in {\mathcal P}(\G_* )$ provide a new
$\mathbb C$-basis for $\SG $.

We define
$\overline{\rho}  \in  \mathcal P ( \G_* )$ by
assigning to $c \in \G_*$ the partition $\rho (c^{-1})$,
which is the composition of $\rho$ with the involution
on $\G_*$ given by $c \mapsto c^{-1}$. It follows from
Proposition~\ref{prop_orth} that

\begin{eqnarray}  \label{eq_inner}
  \langle a_{ - \rho'}', a_{ - \overline{\rho} }' \rangle_{\wt }'
  = \delta_{\rho ', \rho }
    Z_{\rho} \prod_{c \in \G_*} \wt (c)^{l (\rho (c))},
 \quad \rho ', \rho  \in \mathcal P (\G_*).
\end{eqnarray}

\begin{remark}
  $\SGO $ can be characterized as the radical of the
 bilinear form $\langle \ , \ \rangle_{\wt}'$ in $\SG$.
 Thus the bilinear form $\langle \ , \ \rangle_{\wt} '$
 descends to $\SGG$.
\end{remark}
\section{Isometry between $\RG$ and $\SG$}
\label{sect_isom}
\subsection{The characteristic map $\ch$}
  Let $\Psi : \Gn \rightarrow \SG$ be the map defined
by $\Psi (x) = a_{ - \rho}'$ if $x \in \Gn$ is of type $\rho$.
Given $y \in \Gm$ of type $\rho '$, we may regard
$x \times y \in \Gn \times \Gm$ to be in $\G_{n +m}$
of type $\rho \cup \rho '$ so that
$$
  \Psi (x \times y) = \Psi (x) \Psi (y).
$$

      We define a $\mathbb C$-linear map
$ch: \RG \longrightarrow \SG$ by letting
\[
 ch (f )
 = \langle f, \Psi \rangle_{\Gn}
 = \sum_{\rho \in \mathcal P(\G_*)} Z_{\rho}^{-1} f_{\rho}
 a_{ - \rho }',
\]
where $f_{\rho}$ is the value of $f$ at elements of type $\rho$.
The map $ch $ is called the {\em characteristic map}.

We may think of $a_{ -n } (\g ), n >0 , \g \in \G^*$ as
the $n$-th power sum in a sequence of variables
$ y_{ \g } = ( y_{i\g } )_{i \geq 1}$. In this way
we identify the space $\SG$ with the space $\LG$ of
symmetric functions indexed by $ \G^*$ (cf. \cite{M2}).
In particular given a partition $\lambda$
we denote by $s_{\lambda} (\g)$ the Schur
function associated to $y_{\g}$. By abuse of notation,
we denote by $s_{\lambda}(\g)$ the corresponding element
in $\SG$ by the identification of $\SG$ and $\LG$.
For $\lambda \in \mathcal P (\G^*)$, we denote
\begin{eqnarray}  \label{eq_schur}
  s_{\lambda} = \prod_{\g \in \G^*} s_{\lambda(\g)} (\g) \in \SG.
\end{eqnarray}
Then $s_{\lambda}$ is the image
under the isometry $\ch$ of the character of
an irreducible representation $\chi^{\lambda}$ of $\Gn$ (cf. \cite{M2}).

Denote by $c_n (c \in \G_*)$ the conjugacy class in $\Gn$
of elements $(x, s) \in \Gn$ such that $s$ is an $n$-cycle
and the cycle product of $x$ is $ c$.
Denote by $\sigma_n (c )$ the class function on $\Gn$ which takes
value $n \zeta_c$ (i.e. the order of the centralizer of
an element in the class $c_n$) on elements in the class $c_n$
and $0$ elsewhere. For
$\rho = \{ m_r (c) \}_{r \geq 1, c \in \G_*}
\in \mathcal P_n (\G_*)$,
$\sigma_{\rho} = \prod_{r \geq 1, c \in \G_*} \sigma_r (c)^{m_r (c)}$
is the class function on $\Gn$ which takes value
$Z_{\rho}$ on the conjugacy class of type $\rho$ and
$0$ elsewhere. Given $\g \in R(\G)$, we
denote by $\sigma_n (\g )$ the class function on $\Gn $ which takes
value $n \g (c) $ on elements in the class $c_n, c \in \G_*$,
and $0$ elsewhere.

The following lemma is easy to check.
\begin{lemma}  \label{lem_isom}
  The map $ch$ sends $\sigma_{\rho}$ to $a_{ - \rho} '$.
 In particular, it sends $\sigma_n (c)$ to $a_{ -n} (c)$ in $\SG$
 while sending $\sigma_n(\g )$ to $a_{ -n} ( \g )$.
\end{lemma}
\subsection{The characters $\eta_n (\g )$ and $\varepsilon_n (\g )$}

   As we have remarked, the map from $\g  \in \G^*$
to $\eta_n (\g )$ can be extended to be a map
>from $R(\G)$ to $R(\Gn )$. We will need the following
proposition which appeared in \cite{W}
in a more general setting. Eqn.~(\ref{eq_exp}) for $\g \in \G^*$
is well known (cf. \cite{M2}).
We present a complete proof here as we will need
some formulas appearing in the proof later on.

\begin{proposition}  \label{prop_exp}
   For any $\g \in R(\G)$, we have
\begin{eqnarray}
 \sum\limits_{n \ge 0}  \ch ( \eta_n( \g ) ) z^n
  &= & \exp \Biggl( \sum_{ n \ge 1}
      \frac 1n \, a_{-n}(\g )z^n \Biggr), \label{eq_exp} \\
 \sum\limits_{n \ge 0}  \ch ( \varepsilon_n( \g )  ) z^n
  &= & \exp \Biggl( \sum_{ n \ge 1}
      ( -1)^{ n -1} \frac 1n \, a_{-n}(\g )z^n \Biggr). \label{eq_sign}
\end{eqnarray}
\end{proposition}

\begin{proof}
  First assume that $\g \in R(\G)$ is the character of
 a representation $V_{\g }$ of $\G$.
 We calculate the character value of $\eta_n (\g )$
 at $c_n (c \in \G_*)$ first. Take $(g, s) \in \Gn,$ where
 $g =(g_1, \ldots, g_n) \in \G^n$
 and $s$ is an $n$-cycle, say $s = (12 \ldots n)$.
 Denote by $e_1, \ldots, e_k$ a basis of $V_{\g}$, and we write
 $ g e_j = \sum_i a_{ij} (g) e_i $,
 $ a_{ij} (g) \in \mathbb C .$
 It follows that
 \[
  (g, s) (e_{j_1} \otimes \ldots \otimes e_{j_n})
   = g_1 (e_{j_n}) \otimes g_2(e_{j_1})
       \ldots \otimes g_n (e_{j_{n -1}}),
 \]
 in which the coefficient of $e_{j_1} \otimes \ldots \otimes e_{j_n}$
 is
 $$
  a_{j_1 j_n}(g_1) a_{j_2 j_1}(g_2) \ldots a_{j_n j_{n -1}}(g_n).
 $$

 Thus we obtain
 \begin{eqnarray*}
  \eta_n (\g) (c_n)
  & =& \mbox{trace } a(g_n) a(g_{n-1}) \ldots a(g_1)   \\
  & =& \mbox{trace } a(g_n g_{n -1} \ldots g_1) = \g (c)
 \end{eqnarray*}
 since $(g, s) $ is in the conjugacy class $c_n$ which means
 $g_n g_{n -1} \ldots g_1 $ lies in $c \in \G_*$.
 Since the sign character of $\Gn$ takes value
 $(-1)^{n -1}$ at $c_n$, we have
 $$
  \varepsilon_n (\g) (c_n) = ( -1)^{n -1} \g (c).
 $$

 Given $x \times y \in \Gn$, where
 $x \in \G_r$ and $y \in \G_{n -r}$, we clearly have
 $ \eta_n (\g ) (x \times y) = \eta_n (\g) (x) \eta_n (\g)(y).$
 Thus it follows that if $x \in \Gn$ is of type $\rho$, then
 \begin{eqnarray}
  \eta_n (\g ) ( x)
    = \prod_{c\in \G_*} \g (c)^{l (\rho(c))}, \label{eq_term}\\
  \varepsilon_n (\g ) ( x)
    = (-1)^n  \prod_{c\in \G_*} ( - \g (c))^{l (\rho(c))},
   \label{eq_signterm}
 \end{eqnarray}
 where $|| \rho || =n$.
 Putting (\ref{eq_term}) into a generating function, we have
 \begin{eqnarray*}
  \sum\limits_{n \ge 0}  \ch ( \eta_n( \g ) ) z^n
  &= & \sum_{\rho} Z_{\rho}^{ -1}
         \prod_{c\in \G_*} \g (c)^{l (\rho(c))}
          a_{ -\rho (c) }' z^{|| \rho||}                  \\
  &= & \prod_{c\in \G_*} \Bigl ( \sum_{\lambda }
         (\zeta_c^{ -1}\g (c) )^{l (\lambda)}
         z_{\lambda}^{-1} a_{- \lambda} (c) z^{|\lambda|} \Bigr )    \\
  &= & \exp \Biggl  ( \sum\limits_{ n \geq 1}
         \frac1n \sum\limits_{c \in \G_*}
           \zeta_c^{ -1} \g(c) a_{-n} (c) z^n \Biggl )      \\
  &= & \exp \Biggl( \sum_{ n \ge 1}
         \frac 1n \, a_{-n}(\g )z^n \Biggr).
 \end{eqnarray*}
 In a similar manner we can prove (\ref{eq_sign}) by using
 (\ref{eq_signterm}).

 Note that (\ref{eq_sign}) can be obtained from (\ref{eq_exp})
 by substituting $\g$ with $- \g$ and $z$ with $-z$.
 Let  $\beta, \g $ be the characters of two representations
 of $\G$.
 It follows from (\ref{eq_virt}) that

 \begin{eqnarray*}
  & & \sum\limits_{n \ge 0}  \ch ( \eta_n( \beta -\g ) ) z^n  \\
  &=& \Biggl (
       \sum\limits_{n \ge 0}  \ch ( \eta_n( \beta ) ) z^n
      \Biggl )  \cdot
      \Biggl (
       \sum\limits_{n \ge 0} - \ch ( \varepsilon_n( \g ) ) ( -z)^n
      \Biggl )     \\
  &=& \exp \Biggl( \sum_{ n \ge 1}
         \frac 1n \, a_{-n}(\beta ) z^n \Biggr) \cdot
      \exp \Biggl( \sum_{ n \ge 1}
         \frac 1n \, a_{-n}( -\g ) z^n \Biggr)
    \\
  &=& \exp \Biggl( \sum_{ n \ge 1}
         \frac 1n \, a_{-n}(\beta - \g )z^n \Biggr).
 \end{eqnarray*}
Therefore the proposition holds for $\beta - \g$, and so
for any element $\g \in R(\G)$.
\end{proof}

\begin{remark} Formulas (\ref{eq_exp}-\ref{eq_sign}) are equivalent, since
(\ref{eq_sign}) can be obtained (\ref{eq_sign}) from by substituting $\g$
by $-\g$ and $z$ by $-z$.
\end{remark}

\begin{corollary}  \label{cor_char}
   The formula (\ref{eq_term}) holds for any $\g \in R(\G)$.
 In particular $\eta_n (\wt)$ is self-dual.
\end{corollary}

Componentwise, we obtain
\begin{eqnarray*}
  \ch (\eta_n (\g ) )
  &=& \sum\limits_{\la } \frac 1{z_\la }\,
             a_{-\lambda}(\g ), \\
  \ch (\varepsilon_n (\g ))
  &=& \sum\limits_{\la } \frac 1{z_\la }\,
             ( -1)^{ | \la | - l ( \la )} a_{-\lambda}(\g ),
\end{eqnarray*}
where the sum runs over all the partitions $\lambda$ of $n$.
\subsection{Isometry between $\RG$ and $\SG$}
    It is well known that there exists a natural Hopf
algebra structure on the symmetric algebra
$\SG$ with the usual multiplication
and a comultiplication $\Delta$ characterized by
$$
  \Delta ( a_n (\g ))
   = a_n (\g ) \otimes 1 + 1 \otimes a_n (\g ).
$$
Recall that we have defined a Hopf algebra structure
on $\RG$ in Sect.~\ref{sect_wreath}. The following
proposition is easy to check.

\begin{proposition}
  The characteristic map $ \ch: \RG \longrightarrow \SG$
 is an isomorphism of Hopf algebras.
\end{proposition}

    Recall that we have defined a bilinear form
$\langle \  ,  \, \rangle_{\wt }$ on $\RG$ and
a bilinear form on $\SG$ denoted by
$\langle \ , \, \rangle_{\wt }'$. The
lemma below follows from our definition of
$\langle \ , \, \rangle_{\wt }'$ and the
comultiplication $\Delta$.

\begin{lemma}
  The bilinear form $\langle \  ,  \, \rangle_{\wt } '$ on $\SG$
 can be characterized by the following two properties:

 1). $\langle a_{ -n} (\beta ), a_{ -m} (\g ) \rangle_{\wt}'
  = \delta_{n, m} \langle \beta , \g  \rangle_{\wt}' ,
  \quad \beta, \g \in \G^*.$

 2). $ \langle f g , h \rangle_{\wt}'
       = \langle f \otimes g, \Delta h \rangle_{\wt}' ,$
 where $f, g, h \in \SG $, and the bilinear form on the r.h.s
 of 2), which is defined on $\SG \otimes \SG$, is induced
 from $\langle \ , \ \rangle_{\wt}'$ on $\SG$.
\end{lemma}

\begin{theorem}  \label{th_isometry}
  The characteristic map $\ch$ is an isometry from the space
 $ (\RG, \langle \ \ , \ \  \rangle_{\wt } )$ to
 $ (\SG, \langle \ \ , \ \  \rangle_{\wt }' )$.
\end{theorem}

\begin{proof}
  By Corollary~\ref{cor_char}, the character
 value of $\eta_n (\wt )$ at an element $x$ of type $\rho$ is
 $$
  \eta_n (\wt ) ( x)= \prod_{c\in \G_*} \wt (c)^{l (\rho(c))}.
 $$
 Thus it follows by definition of the weighted bilinear form
 \begin{eqnarray*}
   \langle \sigma_{ \rho '}, \sigma_{ \overline{\rho} } \rangle_{\wt }
  & =& \sum_{\mu \in \mathcal P_n (\G_*)}
       Z_{\mu}^{ -1} \wt (c_{\mu }) \sigma_{\rho '} (c_{\mu})
       \sigma_{ \overline{\rho} } (c_{\mu}^{-1} )                      \\
  & =& Z_{\rho}^{ -1}\wt (c_{\rho}) \delta_{\mu, \rho '}
       Z_{\rho '}\delta_{\mu, \rho } Z_{\rho}  \\
  & =& \delta_{\rho, \rho '}
       Z_{\rho}^{ -1}\wt (c_{\rho}) Z_{\rho}Z_{\rho}  \\
  & =& \delta_{\rho, \rho '}
       Z_{\rho} \prod_{c \in \G_*} \wt (c)^{l (\rho (c))}.
 \end{eqnarray*}
 Here $c_{\mu}$ denotes the conjugacy class in $\Gn$ of type $\mu$.
  By Lemma~\ref{lem_isom} and the formula (\ref{eq_inner}), we
 see that
 \[
  \langle \sigma_{ \rho '}, \sigma_{  \overline{\rho} } \rangle_{\wt }
  = \langle a_{- \rho '}, a_{ - \overline{\rho} } \rangle_{\wt } '
  =  \langle \ch (\sigma_{ \rho '}),
             \ch (\sigma_{  \overline{\rho} } ) \rangle_{\wt } '.
 \]
 Since $\sigma_{\rho}, \rho \in \mathcal P(\G_*)$ consist
 a $\mathbb C$-basis of $\RG$, we have established that
 $\ch: \RG \longrightarrow \SG$ is an isometry.
\end{proof}

Thanks to the isometry established above,
we will write $\langle \ , \ \rangle_{\wt}$
for $\langle \ , \ \rangle_{\wt} '$ on $\SG$ from now on.

\begin{remark}
  In the special case when $\wt$ is trivial,
 Theorem~\ref{th_isometry} was established in
 \cite{M1, M2} once we identify $\SG$ with the
 space of symmetric functions parameterized by $\G^*$.
\end{remark}

\begin{remark}
  The standard Hermitian form on $R(\G_n)$ and therefore
on $\RG $ is compatible via the characteristic map $\ch$
with the hermitian form characterized by (\ref{eq_bili})
on $\SG$.

\end{remark}
\section{Vertex operators and $\RG$}
\label{sect_vertex}
\subsection{A $2$-cocycle on an integral lattice}
    Let $Q$ be an integral lattice with a symmetric
bilinear form $\langle \ , \ \rangle$. One can easily
check that the map from $Q$ to $\mathbb Z /2 \mathbb Z$ given by
$\alpha \mapsto ( -1)^{\langle \alpha, \alpha \rangle}$
is a group homomorphism. This homomorphism naturally induces
a $\mathbb Z /2 \mathbb Z$-gradation on $Q$.

Let $\epsilon : Q \times Q \longrightarrow \mathbb C^{\times}$
be such that

\begin{eqnarray}
  \ep (\g , 0 ) & =& \ep (0, \g ) =1 , \label{eq_iden} \\
  \ep (\alpha, \beta ) \ep (\alpha + \beta, \g )
   & =& \ep (\alpha, \beta + \g ) \ep (\beta, \g).  \label{eq_assoc}
\end{eqnarray}
i.e. $\ep$ is a $2$-cocycle of the group $Q$ with values in
$\mathbb C^{\times}$. Another $2$-cocycle $\ep '$ is
equivalent to $\ep$ if and only if there exists
$\epsilon_{\alpha} (\alpha \in \Rz)$ such that
$$
  \ep ' (\alpha, \beta)
  = \ep_{\alpha} \ep_{\beta} \ep_{\alpha +\beta}^{ -1}
    \ep (\alpha, \beta).
$$
The group of equivalent classes of $2$-cocycles is the second
cohomology group $H^2 (Q, \mathbb C^{\times})$.

We introduce
\begin{eqnarray*}
  B_{\epsilon} (\alpha, \beta)
  = \ep (\alpha, \beta) \ep (\beta, \alpha)^{ -1}.
\end{eqnarray*}

Clearly $B_{\ep}$ is skew-symmetric:
$B_{\ep}(\alpha, \beta) = B_{\ep}(\beta, \alpha)^{ -1}.$
It follows from Eqns.~(\ref{eq_iden}) and (\ref{eq_assoc})
that $B_{\ep}$ is bimultiplicative:
\begin{eqnarray*}
  B_{\ep} (\alpha + \beta, \g) &= &
 B_{\ep} (\alpha, \g) B_{\ep} ( \beta, \g), \\
 B_{\ep} (\alpha, \beta+ \g) &= &
 B_{\ep} (\alpha, \beta) B_{\ep} (\alpha, \g).
\end{eqnarray*}
The map $\ep \mapsto B_{\ep}$
gives rise to a homomorphism from $H^2 (Q, \mathbb C^{\times})$
to the group of bimultiplicative skew-symmetric functions on
$Q \times Q$. It is easy to show that
the homomorphism $h$ is indeed an isomorphism.

\begin{proposition}\cite{FK}  \label{prop_cocy}
 There exists a unique $2$-cocyle
 $\ep : Q \times Q \longrightarrow \mathbb C^{\times}$
 up to equivalence such that
 $$
   B_{\ep}(\alpha, \beta) =
  (-1)^{\langle \alpha, \beta \rangle +
        \langle \alpha, \alpha \rangle \langle \beta , \beta \rangle}.
 $$
\end{proposition}
\subsection{Vertex Operators $X ( \g, z)$}
Endowed with the weighted bilinear form
$\langle \ , \ \rangle_{\wt }$, the lattice
$\Rz$ is an integral lattice. We will always associate a $2$-cocycle $\ep$
as in Proposition~\ref{prop_cocy} to the integral lattice
$(\Rz, \langle \ , \ \rangle_{\wt})$ (and its sublattices).
Denote by $\mathbb C[\Rz]$ the group algebra
generated by $e^{\g }$, $\g \in \Rz$.
We introduce two special operators acting on $\mathbb C[ \Rz ]$, namely
a ($\ep$-twisted) multiplication operator $e^{\alpha}$:
 $$
  e^{\alpha }.e^{\beta } = \ep(\alpha, \beta) e^{\alpha +\beta},
 \quad  \alpha, \beta  \in \Rz,
 $$
and a differentiation operator
${\partial_{\gamma }}$:
\begin{eqnarray*}
 {\partial_{\gamma}} e^{ \beta} =
 \langle \gamma, \beta \rangle_{\wt} e^{ \beta},
 \quad  \alpha, \beta  \in \Rz.
\end{eqnarray*}
We extend these two operators acting on
$\mathbb C[\Rz]$ to the following space
$$\FG = \RG \bigotimes \mathbb C[\Rz] ,
$$
by letting them act on the $\RG$ part trivially.

Introduce the operators $ H_{ \pm n}( \g ),
E_{ \pm n} ( \g), \g \in R(\G), n > 0 $
as the following compositions of maps:
\begin{eqnarray*}
  H_{ -n} ( \g ) &:&
    R ( \Gm )
  \stackrel{ \eta_n (\g) \otimes}{\longrightarrow}
    R ( \Gn ) \bigotimes R ( \Gm )
  \stackrel{ {Ind} }{\longrightarrow}
    R ( \G_{n +m} )   \\
  E_{ -n} ( \g ) &:&
    R ( \Gm )
  \stackrel{ \varepsilon_n (\g) \otimes}{\longrightarrow}
    R ( \Gn ) \bigotimes R ( \Gm )
  \stackrel{ {Ind} }{\longrightarrow}
    R ( \G_{n +m} )   \\
  E_n ( \g ) &:&
    R ( \Gm )
   \stackrel{ {Res} }{\longrightarrow}
    R( \Gn) \bigotimes R( \G_{m -n})
   \stackrel{ \langle \varepsilon_n (\g),
         \cdot \rangle_{\wt } }{\longrightarrow}
    R ( \G_{m -n}) \\
  H_n ( \g ) &:&
    R ( \Gm )
   \stackrel{ {Res} }{\longrightarrow}
    R( \Gn) \bigotimes R( \G_{m -n})
   \stackrel{ \langle \eta_n (\g), \cdot \rangle_{\wt} }{\longrightarrow}
    R ( \G_{m -n}) .
\end{eqnarray*}

Define
\begin{eqnarray*}
  H_+ (\g, z) = \sum_{ n > 0} H_{ -n} ( \g ) z^n, &&
  E_+ (\g, z) = \sum_{ n > 0} E_{ -n} ( \g )( -z)^n, \\
  E_- (\g, z) = \sum_{ n > 0} E_{ n} (\g )( -z)^{ -n}, &&
  H_- (\g, z) = \sum_{ n > 0} H_{ n} (\g )z^{ -n} .
\end{eqnarray*}
  Further define operators
$X_n (\g ), n \in {\mathbb Z} + \langle \g, \g \rangle_{ \wt} /2$
by the following generating functions:
\begin{eqnarray}  \label{eq_vo}
 X^+ ( \g, z)
 & \equiv& X ( \g, z) \\
  & =& \sum\limits_{n \in
  {\mathbb Z} + \langle \g, \g \rangle_{ \wt} /2}
 X_n( \gamma) z^{ -n - \langle \g, \g \rangle_{ \wt} /2}  \nonumber     \\
  & =&  H_+ (\g , z) E_- (\g , z) e^{ \g} z^{ \partial_{ \g}}.\nonumber
\end{eqnarray}

Sometimes we denote
\begin{eqnarray*}
 X^- ( \g, z)
  & \equiv & X ( -\g, z)  \\
  & =& \sum\limits_{n \in \mathbb Z
   + \langle \g, \g \rangle_{ \wt}/2}
  X^-_n( \g )z^{-n - \langle \g, \g \rangle_{ \wt}/2}.
\end{eqnarray*}
One easily sees that the operators $X_n (\g )$
are well-defined operators acting on the space $ \FG.$
 We extend the bilinear form
$\langle \ , \ \rangle_{\wt}$ on $\RG$ to $\FG $ by letting
\[
  \langle a .e^{\alpha }, b .e^{ \beta } \rangle_{ \wt} =
   \langle a , b \rangle_{\wt} \delta_{ \alpha, \beta }, \quad
    a, b \in \RG, \, \alpha, \beta \in \Rz.
\]
We extend the $\mathbb Z_+$-gradation
on $\RG$ to a $\frac12 \mathbb Z_+$-gradation on
$\FG$ by letting
\begin{eqnarray*}
  \deg  a_{ -n} (\g ) = n , \quad
  \deg e^{\g } = \frac12 \langle \g , \g \rangle_{\wt}.
\end{eqnarray*}

Similarly we extend the bilinear form
$\langle \  , \ \rangle_{\wt }$  to the space
$$
  \VG = \SG \bigotimes \mathbb C [ \Rz]
$$
and extend the $\mathbb Z_+$-gradation on $\SG$
to a $\frac12 \mathbb Z_+$-gradation on $\VG$.

We extend the characteristic map $\ch$ to an isometry
>from $\FG$ to $\VG$ by
the identity operator on the factor $\mathbb C [ \Rz]$,
which is denoted again by $\ch$.
\subsection{Heisenberg algebra and $\RG$}
    We define $ \widetilde{a}_{ -n} (\gamma), n >0$ to be a map
>from $\RG$ to itself by the following composition
\[  R (\Gm) \stackrel{ \sigma_n ( \g ) \otimes }{\longrightarrow}
  R(\Gn) \bigotimes R (\Gm)  \stackrel{{Ind} }{\longrightarrow}
  R ( {\Gamma}_{n +m}).
\]
We also define $ \widetilde{a}_{ n} (\gamma), n >0$ to be a map from $\RG$
to itself
as the composition
\[
  R (\Gm)  \stackrel{ Res }{\longrightarrow}
   R(\Gn)\bigotimes R ( {\G }_{m -n})
 \stackrel{ \langle \sigma_n ( \g), \cdot \rangle_{\wt}}{\longrightarrow}
 R ( {\G }_{m -n}).
\]

 We denote by $\RGO $ the radical of the bilinear form
$\langle \ , \ \rangle_{\wt}$ in $\RG$ and denote by
$\RGG$ the quotient $\RG / \RGO$. $\RGG$ inherits the
bilinear form $\langle \ , \ \rangle_{\wt}$ from $\RG$.

\begin{theorem}  \label{th_heis}
  $\RG$ is a representation of the Heisenberg algebra
 $\hg $  by letting $ a_n (\g )$
 $( n \in \mathbb Z \backslash 0)$ act as $ \widetilde{a}_{ n} (\g )$,
 $a_0 (\g )$ as $0$ and $C$ as $1$. $\RGO$ is a subrepresentation
 of $\RG$ over $\hg$ and the quotient $\RGG $
 is irreducible. The characteristic map
 $\ch$ is an isomorphism of $\RG$ (resp. $\RGO$, $\RGG$)
 and $\SG$ (resp. $\SGO$, $\SGG$) as representations over $\hg$.
\end{theorem}

\begin{remark}
 For a $\G$-space $X$, the wreath product $\Gn$ acts on
 the $n$-th direct product $X^n$. A Heisenberg algebra was
 defined in \cite{W} to act on a direct sum of equivariant K-theory
 $\bigoplus_{n \geq 0} K_{\Gn } (X^n) \bigotimes \mathbb C$.
 In view of Corollary~\ref{cor_char},
 our Heisenberg algebra here is a special case of
 the Heisenberg algebra constructed in \cite{W}
 when $X$ is a point, cf. \cite{W} for a proof.
\end{remark}
\subsection{Vertex operators and Heisenberg algebra $\hg $}
    The relation between the vertex operators defined
in (\ref{eq_vo}) and the Heisenberg algebra $\hg $ is revealed
in the following theorem.

\begin{theorem} For any $\g \in R(\G)$, we have
  \begin{eqnarray*}
   \ch \bigl ( H_+ (\g, z) \bigl )
   &=& \exp \biggl ( \sum\limits_{ n \ge 1} \frac 1n \,
    a_{-n} ( \g ) z^n \biggr ), \\
   \ch \bigl ( E_+ (\g, z) \bigl )
   &=& \exp\biggl (- \sum\limits_{n\ge 1} \frac 1n  \,
    a_{-n}(\gamma)z^n\biggr ),  \\
   \ch \bigl ( H_- (\g , z) \bigl )
   &=& \exp \biggl ( \sum\limits_{n \ge 1}\frac 1n \,
    a_n (\g ) z^{-n}\biggr ),         \\
   \ch \bigl ( E_- (\g , z) \bigl )
   &=& \exp\,\, \biggl ( -\sum\limits_{ n \ge 1}
    \frac 1n \,{ a_n( \g )} z^{ -n} \biggr ) .
  \end{eqnarray*}
\end{theorem}

\begin{proof}
  The first and second identities were essentially
 established in Proposition~\ref{prop_exp} together with
 Lemma~\ref{lem_isom}. The only minor difference is that
 the components appearing in these two identities are
 regarded as operators acting on $\RG$ and $\SG$.

 Note that by definition
 the adjoints of $E_+ (\g, z)$ and $H_+ (\g , z)$
 with respect to the bilinear form
 $\langle \ , \ \rangle_{\wt}$
 are $E_- (\g , z^{ -1})$ and $ H_- (\g , z^{ -1})$, respectively.
 The third and fourth identities are obtained by applying
 the adjoint functor to the first two identities.
\end{proof}

Putting all pieces together, we have
\begin{eqnarray*}
  && \ch \bigl ( X( \g , z)\bigl )  \\
  &= & \exp \biggl ( \sum\limits_{ n \ge 1}
  \frac 1n \, a_{-n} ( \g ) z^n \biggr ) \,
  \exp \biggl ( -\sum\limits_{ n \ge 1}
  \frac 1n \,{ a_n( \g)} z^{ -n} \biggr )
  e^{ \g} z^{ \partial_{\g }}.
\end{eqnarray*}
\section{Vertex representations and the McKay correspondence}
\label{sect_ade}
\subsection{Product of two vertex operators}
        We first note that
$$
  X_n(\gamma)^* = X_{-n}(- \gamma),
 \quad n\in \mathbb Z +\langle \g, \g \rangle_{\wt}/2.
$$
The normal ordered product $: X(\alpha, z) X(\beta, w) :$,
$\alpha, \beta \in R(\G)$ of two vertex operators is defined
as follows:
\begin{eqnarray*}
 &  & : X(\alpha, z) X(\beta, w) :  \\
 & =&
  \exp
  \biggl ( \sum\limits_{ n \ge 1}
   \frac 1n \,
    \bigl(
      a_{-n} ( \alpha ) z^n + a_{-n} ( \beta ) w^n
    \bigr)
  \biggr ) \cdot                              \\
  & &
 \cdot
 \exp
  \biggl ( -\sum\limits_{ n \ge 1}
  \frac 1n \,
    \bigl(
      a_n( \g) z^{ -n}+ a_n ( \beta ) w^{ -n}
    \bigr)
  \biggr )
  e^{ \alpha +\beta} z^{ \partial_{\alpha }} w^{ \partial_{\beta}}.
\end{eqnarray*}
 The term $(z -w)^{ \langle \alpha, \beta \rangle_{\wt}}$
in the following theorem is understood as the formal series expansion
(cf. \cite{FLM})
$$
  z^{ \langle \alpha, \beta \rangle_{\wt}}
\sum_{ k \geq 0} \binom{\langle \alpha, \beta \rangle_{\wt}}{k}
(- z^{-1}w)^k.
$$
This applies to similar expressions in similar contexts below.
\begin{theorem}  \label{th_ope}
 For $\alpha, \beta \in \Rz$, we have the following identity
 for a product of two vertex operators:
 \begin{eqnarray*}
  X(\alpha, z) X(\beta, w) & =& \ep (\alpha, \beta)
  :X(\alpha, z) X(\beta, w):
  (z -w)^{ \langle \alpha, \beta \rangle_{\wt}}.
 \end{eqnarray*}
         %
         %
\end{theorem}

\begin{proof}
   We first compute
 \begin{eqnarray*}
   && E_- (\alpha, z) H_+ (\beta, w)    \nonumber \\
   & =& \exp\,\biggl ( -\sum\limits_{ n \ge 1}
    \frac 1n \,{ a_n(\alpha)} z^{ -n} \biggr )
   \exp \biggl ( \sum\limits_{ n \ge 1} \frac 1n \,
    a_{-n} (\beta) w^n \biggr )    \nonumber  \\
  & =&  \exp \biggl ( \sum\limits_{ n \ge 1} \frac 1n \,
       a_{-n} (\beta) w^n \biggr )
       \exp\,\biggl  ( -\sum\limits_{ n \ge 1}
       \frac 1n \,{ a_n(\alpha)} z^{ -n} \biggr ) \cdot   \\
  &  & \cdot \exp \, \biggl  ( - \langle \alpha, \beta \rangle_{\wt}
       \sum\limits_{ n \ge 1} \frac 1n \,
       z^{ -n} w^n \biggr )  \nonumber  \\
  & =&  \exp \biggl ( \sum\limits_{ n \ge 1} \frac 1n \,
        a_{-n} (\beta) w^n \biggr )
        \exp\,\biggl  ( -\sum\limits_{ n \ge 1}
        \frac 1n \,{ a_n(\alpha)} z^{ -n} \biggr )
        (1 - z^{ -1} w )^{\langle \alpha, \beta \rangle_{\wt} } .
      \nonumber
  \end{eqnarray*}
 The second identity above uses the formula
 $$
  \exp (A) \exp (B) = \exp (B) \exp (A) \exp ( [A, B])
 $$
 since $[A, B]$ here commutes with $A$ and $B$.

 On the other hand, we have
 \[
  z^{\partial_{\g_i}} e^{ \g_j}
  = z^{a_{ij}} e^{\g_j} z^{\partial_{\g_i}} .
 \]
  Combining with Eqns.~(\ref{eq_iden}) and (\ref{eq_assoc}) we compute
  \begin{eqnarray*}
   &&  X(\alpha, z) X(\beta , w)   \\
  & =& H_+ (\alpha, z) E_- (\alpha, z) e^{\alpha} z^{\partial_{ \alpha_i}}
       H_+ (\beta, w) E_- (\beta, w) e^{\beta} w^{ \partial_{\beta}}  \\
  & =& \ep (\alpha, \beta)
       H_+ (\alpha, z) H_+ (\beta, w)
       E_- (\alpha, z) E_- (\beta, w) \cdot \\
  &  & \cdot
       e^{\alpha + \beta} z^{\langle \alpha, \beta \rangle_{\wt}}
       z^{ \partial_{\alpha}} w^{ \partial_{\beta}}
       \Bigl (1 - \frac wz \Bigr )^{\langle \alpha, \beta \rangle_{\wt}}  \\
  & =& \ep (\alpha, \beta)
       :X(\alpha, z) X(\beta, w):
        (z -w)^{\langle \alpha, \beta \rangle_{\wt}} .
  \end{eqnarray*}
\end{proof}

The following proposition is easy to check.

\begin{proposition}  \label{prop_prim}
  Given $\alpha \in R(\G), \beta \in \Rz$, we have
  $$
   [ a_n (\alpha ), X(\beta, z)]
   = \langle \alpha, \beta \rangle_{\wt} X(\beta, z) z^n.
  $$
\end{proposition}
\subsection{Affine Lie algebra $\loopg$ and toroidal Lie algebra $\hhg$}
  Let $\mathfrak g$ be a complex simple Lie algebra of ADE type
of rank $r$,
with a root system $\overline{\Delta}$ and a set of simple roots
$\alpha_1. \ldots, \alpha_r$.
Denote by $\theta$ the highest root.
Denote by $e_{\alpha}, f_{\alpha}, h_{\alpha} $,
$\alpha \in \overline{\Delta}$ the Chevalley generators.
Sometimes we will write $e_i, f_i$ and $h_i$
for $e_{\alpha_i}, f_{\alpha_i}$ and $ h_{\alpha_i}$ respectively.
Denote by $( \ , \ )$ the invariant bilinear
form on $\mathfrak g$ normalized by letting $ (\theta, \theta) = 2$.
Let $\mathbb C [t, t^{-1}]$
be the space of Laurent polynomials in an indeterminate $t$.
The affine algebra $\loopg$ is the universal
central extension of
$\mathfrak g \bigotimes \mathbb C [t, t^{-1}]$
by a one-dimensional center with a generator $C$.
As a vector space,
$\widehat{\mathfrak g} =
 \mathfrak g \bigotimes \mathbb C [t, t^{-1}]
  \bigoplus \mathbb C C.
$
We denote by
$$
 a (n) = a \otimes t^n, \quad a \in \mathfrak g, n \in \mathbb Z.
$$
The commutation relations in $\loopg$ are given as follows:
\begin{eqnarray*}
  [ a(n), b(m)] & =& [a, b] (n+m) +n \delta_{n, -m} (a, b) C , \\
  {[C, a(n)]}   & =& 0, \quad a, b \in \mathfrak g, n, m \in \mathbb Z.
\end{eqnarray*}

We denote by $A = (a_{ij})_{ 0 \leq i,j \leq r}$ the
affine Cartan matrix associated to $\loopg$. The matrix
$(a_{ij})_{ 1 \leq i,j \leq r}$ obtained by deleting the first
row and column of $A$ is the Cartan matrix of $\mathfrak g$.
We fix a $2$-cocycle $\varepsilon$ on the affine root lattice
given in Proposition~\ref{prop_cocy}.

The so-called {\em basic} representation  $V$ of $\loopg$ is
an irreducible highest weight representation generated by
a highest weight vector which is annihilated by
$ a(n), n \geq 0, a \in \mathfrak g$; $C$ acts on $V$ as the
identity operator.

We denote by $\hhg$ the toroidal Lie algebra over $\mathbb C$
(associated to $\mathfrak g$) with the
following presentation \cite{MRY}:
generators are
$$
 C, h_i (n), x_n (\pm \alpha_i), n \in \mathbb Z, i =0, \ldots, r;
$$
relations are given by:
$C$ is central, and
\begin{eqnarray}
  [ h_i (n), h_j (m)]
  & =& n a_{ij} \delta_{n, -m}C,    \nonumber  \\
  {[h_i (n), x_m (\pm \alpha_j)]}
  & =& \pm a_{ij} x_{n +m} (\pm \alpha_j),  \nonumber    \\
  {[x_n (\alpha_i), x_m (- \alpha_j)]}
  & =& \delta_{ij} \ep(\alpha_i, -\alpha_i)
    \{ h_i (n +m) + n \delta_{n, -m} C \},
      \label{eq_pres}  \\
  {[x_n (\pm \alpha_i), x_m (\pm \alpha_i)]  }
  & =& 0 ,      \nonumber   \\
  (ad \, x_0 (\pm \alpha_i))^{1 -a_{ij}} x_m (\pm \alpha_j)
  & =& 0,  \quad (i \neq j), \nonumber
\end{eqnarray}
where $n, m \in \mathbb Z$, $i, j = 0, 1, \ldots, r$.

There exists a surjective homomorphism from the toroidal algebra $\hhg$ to
the double loop algebra $\mathfrak g \bigotimes  \mathbb C
[ s, s^{-1}, t, t^{-1}]$
given by:
\begin{eqnarray*}
  C \mapsto 0,  &&
  h_i (k) \mapsto h_i \otimes s^k, \quad  i = 0, \ldots, r; \\
  x_n (\alpha_i) \mapsto e_i \otimes s^n,   &&
  x_n (- \alpha_i) \mapsto f_i \otimes s^n, \quad i = 1, \ldots, r; \\
  x_n (\alpha_0) \mapsto f_{\theta} \otimes s^n t, &&
  x_n (- \alpha_0) \mapsto e_{\theta} \otimes s^n t^{-1}, \quad n \in \mathbb Z.
\end{eqnarray*}
\subsection{A new form of McKay correspondence}
  In this subsection we set $\G$ to be a finite
subgroup of $SU_2$ and $\wt$ to be the
virtual character $ 2 \g_0 - \pi$ of $\G$, where
$\pi$ is the character of the two-dimensional natural representation
in which $\G$ embeds in $SU_2$.
We recall that the matrix $A = (a_{ij})_{0 \leq i,j \leq r}$
in Sect.~\ref{subsec_mcka}
is the Cartan matrix for the corresponding affine Lie
algebra $\loopg$.

The following theorem provides a finite group
realization of the vertex representation of
the toroidal Lie algebra $\hhg$ on $\FG\cong \VG$.

\begin{theorem}  \label{th_mck}
 A vertex representation of the toroidal Lie algebra $\hhg$
 is defined on the space $\FG$ by letting
 \begin{eqnarray*}
  x_n (\alpha_i ) \mapsto X_n (\g_i),  &&
  x_n (- \alpha_i ) \mapsto X_n (- \g_i),\\
  h_i (n) \mapsto a_n (\g_i), && C \mapsto 1,
 \end{eqnarray*}
 where $n \in \mathbb Z, 0 \leq i \leq r$.
\end{theorem}

\begin{proof}
 The identity for the product
 of two vertex operators in Theorem~\ref{th_ope}
 implies the commutation relations between
 the components of vertex operators either by
 formal calculus or the Cauchy residue formula
 (see \cite{FLM}). Taking into account
 Proposition~\ref{prop_prim} one deduces that the
 correspondence defined in Theorem~\ref{th_mck}
 indeed realizes a representation of the toroidal algebra $\hhg$,
 (cf. \cite{F2, MRY}).
\end{proof}

 Recall that $\delta = \sum_{i =0}^r d_i \g_i$ generates the
one-dimensional radical $R^0_{\mathbb Z}$ of the bilinear form
$\langle \ , \ \rangle_{\wt}$ in $\Rz$. The lattice $\Rz$ in this case
can be identified with the root lattice
for the corresponding affine Lie algebra.
The quotient lattice $\Rz / R^0_{\mathbb Z}$ inherits
a positive definite integral bilinear form.
Denote by $\Gbar$ the set of non-trivial irreducible characters of $\G$:
$$
  \Gbar = \{ \g_1, \g_2, \ldots, \g_r \}.
$$
Denote by $\Rzz$ the sublattice in $\Rz$
generated by $\Gbar$. The subset of elements in $\Rzz$
of length $\sqrt{2}$ can be identified with the finite
root system $\overline{\Delta}$.
Equipped with bilinear form $\langle \ , \ \rangle_{\wt}$,
$\Rzz $ is obviously
isomorphic to $\Rz / R^0_{\mathbb Z}$, which is in turn
identified with the root lattice of Lie algebra $\mathfrak g$.
Furthermore $\Rz$ can be written as a direct sum of two
integral lattices $ \Rzz \bigoplus \mathbb Z \delta$.

Denote by $Sym (\Gbar)$ the symmetric algebra generated by
$a_{-n} (\g_i),$ $ n >0, $ $i =1, \ldots , r$. Equipped with the bilinear
form $\langle \ , \ \rangle_{\wt}$, $Sym (\Gbar)$
is isometric to $\SGG$ which is in turn isometric
to $\RGG$ as well.

We define
\begin{align*}
 \VGG  &= \SGG \bigotimes \mathbb C [ \Rz / R^0_{\mathbb Z}]
 \cong Sym (\Gbar) \bigotimes \mathbb C [ \Rzz],\\
\FGG  &= \RGG \bigotimes \mathbb C [ \Rzz]
 \cong \VGG.
\end{align*}
 The space $\VG$ associated to the lattice
$\Rz$ is isomorphic to the tensor
product of the space $\VGG$
associated to $\Rzz$ and the space associated to
the rank $1$ lattice $\mathbb Z \delta$
equipped with the zero bilinear form.

The identity for a product of vertex operators $X (\g, z)$
associated to $\g \in \overline{\Delta}$
(cf.~Theorem~\ref{th_mck}) imply that
$\VGG$ affords a realization of the vertex representation
of $\loopg$ on $\VGG$ \cite{FK, S1}.
The following theorem provides
a direct link from the finite group $\G \in SU_2$
to the affine Lie algebra $\loopg$ and thus
it can be regarded as a new form
of McKay correspondence.

\begin{theorem}
 The operators $X_n(\g ), \g \in \overline{\Delta}, a_n (\g_i),
 i =1, 2, \ldots, r, n \in \mathbb Z$
 define an irreducible representation of the affine Lie algebra $\loopg$
 on $\FGG$ isomorphic to the basic representation. In particular
 the operators $X_0(\g ), \g \in \overline{\Delta}, a_0 (\g_i),
 i =1, 2, \ldots, r$ define a representation of
 the simple Lie algebra $\mathfrak g$.
 $\VGG$ carries a distinguished basis
 $\chi^{\rho} \otimes e^{\g}, \rho \in \mathcal P (\Gbar),
  \g \in \Rzz.$
\end{theorem}
\section{Vertex operators and irreducible characters of $\Gn $}
\label{sect_char}
     In this section we specialize $\wt$ to be
the trivial character $\g_0$ of $\G$.
\subsection{Algebra of vertex operators for $\wt = \g_0$}
      In this case the weighted bilinear form reduces
to the standard one $\langle \ , \ \rangle$
and $\Rz$ is isomorphic to the lattice $\mathbb Z^n$
with the standard integral bilinear form. Recall that
$\langle \g_i , \g_j \rangle = \delta_{ij}$. The $2$-cocycle
$\ep$ can be chosen by letting
$ \ep (\g_i, \g_j ) = 1$ if $i \leq j$ and $-1$ if $i >j$.
It follows by standard argument that
the identity for product of vertex operators given in
Theorem~\ref{th_ope} implies the following
anti-commutation relations (cf. \cite{F1}).
The bracket
$\{ \ , \ \}$ below denotes the anti-commutator.

\begin{theorem}  \label{th_cliff}
  The operators $X^+_n (\g_i), X^-_n (\g_i)$
 $(n \in \mathbb Z + \frac12, 0 \leq i \leq r)$
 generate a Clifford algebra:
 \begin{eqnarray}
  \{X^+_m( \g_i), X^+_n(\g_j) \} &= & 0, \nonumber   \\
  \{ X_m^-( \g_i), X_n^-(\g_j ) \} &= & 0,  \nonumber  \\
  \{X^+_m(\g_i), X^-_n(\g_i)\} &= & \delta_{ij}\delta_{m, -n}.
  \label{eq_delta}
 \end{eqnarray}
\end{theorem}
\subsection{Character tables of $\Gn$ and vertex operators}
We now construct a special orthonormal basis in $\VG$ and then
interpret them as irreducible characters of $\Gn$.
First we have the following simple lemma, cf. \cite{J3}
for a proof.

\begin{lemma}
  \label{L:lattice}
 For $m\in \mathbb Z + 1/2$, $\alpha \in R( \G )$ and
 $ \g \in \G^*$, we have
 \begin{eqnarray*}
  X_{m}(\g ) e^{\alpha}
   =& \delta_{m, -\langle \alpha , \g \rangle -1/2}e^{ \a + \g },&
    m \geq - \langle \g, \a  \rangle - 1/2,            \\
  X_{m}^-(\g ) e^{\alpha}
   =& \delta_{m, \langle \a, \g \rangle - 1/2} e^{\a -\g }, &
    m \geq \langle \g, \a  \rangle - 1/2 .
 \end{eqnarray*}
\end{lemma}

Consequently we have
\begin{align*}
 X_{-m -\langle \g, \a  \rangle + 1/2}(\gamma)
 X_{-m -\langle \g, \a  \rangle + 3/2}(\gamma)
  \cdots
 X_{-\langle \g, \a \rangle - 1/2}(\g )e^{\alpha}
 &= e^{ \a + m \g } ,                           \\
 X^-_{-m +\langle \g, \a \rangle + 1/2}(\gamma)
 X^-_{-m +\langle \g, \a \rangle + 3/2}(\gamma)
  \cdots
 X^-_{\langle \g, \a \rangle - 1/2}(\g )e^{\alpha}
 &= e^{ \a - m \g } .
\end{align*}

For a $m$-tuple index $\phi=(\phi_1, \cdots, \phi_m)\in
(\mathbb Z+1/2)^m$ we denote
\begin{equation*}
X_{\phi}(\g )=X_{\phi_1}(\g )\cdots X_{\phi_m}(\g )
\end{equation*}
We will use similar notations for other vertex operators.

Let $\delta_l=(l-1, l-2, \cdots, 1, 0)$ be the special partition
of length $l$. We often omit the subscript if the meaning is clear from
the context. Given $ \la \in \mathcal P(\G^*)$, we define
$$
  \omega(\la) = \sum_{\g \in \G^*}l(\la(\g )) \g \in \Rz ,
$$
and
$$
 s_{\lambda, \alpha}
  = \prod_{\g \in \mathcal \G^*} X_{-\la( \g ) -\delta
   -( \langle \g , \a \rangle +1/2){\bf 1}} (\g).e^{\a},
$$
where the order of $\prod_{\g \in \mathcal \G^*}$ is understood as
a product over $\g_0, \g_1, \ldots, \g_r$ from the left to the right.

Recall that $s_{\lambda} \in \SG$ for $ \lambda \in \mathcal P ( \G^*)$
is defined in Eq.~(\ref{eq_schur}).

\begin{theorem}
 The vectors $s_{\lambda, \alpha}$ for
 $ \lambda=(\la(\g ))_{\g \in \G^*} \in \mathcal P (\G^*),
 \alpha \in \Rz$
 form an orthonormal basis in the vector space $\mathcal F_{\G}$, where
 ${ \bf 1}=(1, \cdots, 1)$ and $\delta $ are both of length
 $l =l(\la(\g ))$. The isometry $\ch$ maps the basis vector
 $s_{\lambda, \alpha}$ to
 $
 s_{\lambda} e^{\a + \omega(\la)} \in \VG.$
\end{theorem}

\begin{proof}
  The Clifford algebra structure (\ref{eq_delta})
 implies that the nonzero elements
 $$
  \prod_{\gamma\in \G^*}X_{-n_1}(\gamma)\cdots
  X_{-n_l}(\gamma)e^{\alpha}
 $$
 of distinct indices generate a spanning set in the space $\VG $.
 To see that they correspond to Schur functions we compute
 in a way similar to Theorem~\ref{th_ope} that
 \begin{align*}
   &X (\g , z_1)\cdots X (\g , z_l)e^{\a}  \\
   & = :X (\g , z_1) \cdots X (\g , z_l):
      \prod_{i <j} (1 - z_j {z_i}^{ -1})
      z^{\delta + ( 1/2 + \langle \g , \a \rangle ){\bf 1}}e^{\a} .
 \end{align*}
 The method in \cite{J1, J2, J3} implies that this is
 exactly the generating function of Schur functions
 under the isomorphism $\ch$. Since the vertex operators
 corresponding to different $\g \in \G^*$
 commute up to a sign it follows
 that the isometry $\ch$ maps the basis vector
 $s_{\lambda, \alpha}$ to
 $ s_{\lambda} e^{\a +\omega (\la )} \in \VG.$

 The orthonormality follows readily from the Clifford
 algebra commutation relations in Theorem~\ref{th_cliff}.
\end{proof}

It is a routine computation either by
symmetric functions or vertex operator
calculus (cf. \cite{J1}) that the complete character
table of $\Gn$ is given by matrix coefficients of
the vertex operators $X(\g, z)$.

\begin{theorem}\label{T:main}
 Given partition-valued functions $\lambda \in \mathcal P ( \G^*)$
 and $\mu \in \mathcal P (\G_*)$,
 the matrix coefficient
 \begin{eqnarray*}
  \langle
    s_{\lambda, -\omega(\la) }, a_{-\mu} '
  \rangle
 \end{eqnarray*}
 is equal to the character value of the irreducible character
 of $\chi^{\la}$ associated to $\lambda$
 at the conjugacy class of type $\mu$ in $\Gn$.
\end{theorem}

\bibliographystyle{amsalpha}

\end{document}